\documentclass[12pt,english]{article}
\usepackage{theorem}
\usepackage{babel}
\usepackage{a4wide}
\usepackage{graphics}
\usepackage{epsfig}
\usepackage{amssymb}
\usepackage{latexsym}
\usepackage{amsmath}

\newtheorem{thh}{Theorem}[section]
\newtheorem{df}[thh]{Definition}

\newtheorem{lem}[thh]{Lemma}

\newtheorem{ques}[thh]{Question}

\title{Group actions, $k$-derivations and finite morphisms}
\author{Philippe Bonnet}
\date{}

\newcommand{\dem}{{\em Proof: }}
\newcommand{\qed}{\begin{flushright} $\blacksquare$\end{flushright}}

\newcommand{\OX}{{\cal{O}}_X}
\newcommand{\XG}{{X/\!/G}}
\newcommand{\OY}{{\cal{O}}_Y}
\newcommand{\CC}{\mathbb{C}}

\newcommand{\FF}{{\cal{F}}}

\begin{document}
\maketitle

\begin{center}
Mathematisches Institut, Universit\"at Basel \\
Rheinsprung 21, 4051 Basel, Switzerland \\ e-mail:
Philippe.bonnet@unibas.ch
\end{center}

\begin{abstract}
Let $G$ be an affine algebraic group over an algebraically closed field $k$ of characteristic zero. In this paper, we consider finite
$G$-equivariant morphisms $F:X\rightarrow Y$ of irreducible affine varieties. First we determine under which conditions on $Y$ the
induced map $F^G:X/\!/G\rightarrow Y/\!/G$ of quotient varieties is also finite. This result is reformulated in terms of kernels
of derivations on $k$-algebras $A\subset B$ such that $B$ is integral over $A$. Second we construct explicitly two examples of
finite $G$-equivariant maps $F$. In the first one, $F^G$ is quasifinite but not finite. In the second one, $F^G$ is not even quasifinite.
\end{abstract}

\section{Introduction}

We begin by recalling a few definitions from Geometric Invariant Theory; for more details, see \cite{B-S}.
Let $G$ be a connected affine algebraic group over an algebraically closed field $k$ of characteristic zero. An affine
$G$-variety $X$ is
an affine variety together
with an algebraic action of $G$, i.e. a regular map $\varphi: G\times X\rightarrow X,\;(g,x)\mapsto g.x$ that
also defines an action of $G$
on $X$. A $G$-equivariant morphism $F:X\rightarrow Y$ of affine $G$-varieties is a regular map
such that $F(g.x)=g.F(x)$ for all $(g,x)\in G\times X$.
For any affine $G$-variety $X$, denote by ${\cal{O}}^G _X$ the ring of
invariants of $G$. The algebraic quotient $X/\!/G$ is defined as
the scheme:
$$
X/\!/G=Spec({\cal{O}}^G _X)
$$
This notion is directly related to Hilbert's fourteenth Problem in Geometric Invariant Theory (see \cite{Van}).
It is well-known that $\XG$ is affine if $G$ is linearly reductive (see \cite{B-S}) but that it neednot be affine in
general (see for instance \cite{Na} or \cite{D-F}). Let
$F:X\rightarrow Y$ be a $G$-equivariant morphism of
affine $G$-varieties. Since $F^{*}({\cal{O}}^G _Y) \subset {\cal{O}}^G _X$,
the map $F$ induces a natural morphism $F^G: X/\!/G \rightarrow Y/\!/G$. In this paper, we are interested in the
following question:

\begin{ques} \label{one}
If $F$ is a finite morphism, under which conditions is $F^G$ also finite?
\end{ques}
In \cite{Van}, p. 227, Van den Essen gave an example of a $G_a$-equivariant finite morphism $F: k^3 \rightarrow C\times k^2$, where
$C$ is a cuspidal curve. In this example, the quotient $k^3/\!/G_a$ is affine but $C\times k^2/\!/G_a$ is not, from which follows
that $F^{G_a}$ is not finite. Note that, in the
case of a unipotent group $G$, the map $F^G$ is always quasifinite (see lemma \ref{uni})

In any case, the Lie algebra ${\mathfrak{g}}$ of $G$ acts like a
collection of $k$-derivations on $\OX$. Since $k$ has characteristic zero, the ring ${\cal{O}}^G _X$ coincides
with the kernel of this collection (see \cite{Kr}). This invites us to reformulate question \ref{one} in a more
general setting. Consider two integral $k$-algebras $A,B$ of finite type,
where $char(k)=0$ and $A\subset B$. Let ${\cal{F}}$ be any family of $k$-derivations on $B$ which preserve $A$, i.e.
$d(A)\subset A$ for any $d\in \FF$. Denote by:
$$
A^{\FF}=\cap_{d\in \FF} \ker d_{|A} \quad \quad \quad \quad B^{\FF}=\cap_{d\in \FF} \ker d_{|B}
$$
\begin{ques} \label{first}
If $B$ is integral over $A$, under which conditions is $B^{\FF}$ is integral over $A^{\FF}$?
\end{ques}
Note that, if $A,B$ are fields instead of finitely generated $k$-algebras, then $B^{\FF}$ is obviously algebraic, hence integral
over $A^{\FF}$. This can be easily proved by using Differential Galois Theory (see \cite{Ko}).

In this paper, we are going to give an answer to these questions, in terms of the properties of the scheme $Spec(A)$,
and then produce some new counterexamples. We begin with some definitions. Let $A$ be a finitely generated $k$-algebra.
A closed point $x$ of $Spec(A)$ (considered as a maximal ideal of $A$) is {\em a non normal singular point} if the
local ring $A_x$ is not integrally closed in its fraction field. In this case, $Spec(A)$ is not smooth at $x$. However,
any singularity does not have to be non normal. For instance, the surface $S$ in $\CC^3$ given by the equation
$xz-y^2=0$ has a unique normal singularity at the origin.

\begin{df}
A finitely generated $k$-algebra $A$ has isolated non normal singular points if the set of non normal singular points of $Spec(A)$ is finite.
\end{df}

\begin{thh} \label{premier}
Let $A\subset B$ be two finitely generated integral $k$-algebras, where $char(k)=0$. Let
$\FF$ be a family of $k$-derivations on $B$ such that $d(A)\subset A$ for any $d\in \FF$. If $B$ is integral over $A$
and if $A$ has isolated non normal singular points, then $B^{\FF}$ is integral over $A^{\FF}$.
\end{thh}
This result is the best one can expect, regarding the singularities of $Spec(A)$. More precisely, there exist examples of actions on
varieties with as few singularities as possible, and for which the conclusion of theorem \ref{premier} fails. We construct such examples
as follows. For any positive integers $n,m$, set $B_{n,m}=k[x_1,..,x_n,y_1,..,y_m,z]$, and define $A_{n,m}$ as the $k$-subalgebra of
$B_{n,m}$ generated by:
\begin{itemize}
\item{the monomials $y_1,...,y_m,z$,}
\item{the polynomials $x_i ^2 + x_i z$ and $x_i ^3 + x_i ^2 z$, for any $i=1,...,n$,}
\item{the monomials $x_1 ^{i_1}...x_n ^{i_n}y_j$, where $j=1,...,m$ and $i_k\leq 1$ for any $k=1,...,n$.}
\end{itemize}
Note that $B_{n,m}$ is integral over $A_{n,m}$. Indeed,
the monomials $y_1,..,y_m,z$ belong to $A_{n,m}$. If $t_i=x_i ^2 + x_i z$, then
$t_i$ belongs to $A_{n,m}$ and each $x_i$ satisfies the relation $x_i ^2 + zx_i -t_i=0$. Since $B_{n,m}$ is normal and has the
same fraction field as $A_{n,m}$, $B_{n,m}$ is the integral closure of $A_{n,m}$.
We set $X_{n,m}=k^{n+m+1}$ and $Y_{n,m}=Spec(A_{n,m})$. The inclusion $A_{n,m}\subset
B_{n,m}$ induces the so-called normalization morphism:
$$
F: X_{n,m} \longrightarrow Y_{n,m}
$$
which is finite. Let $G_1$ be the additive group $G_a(k)$ and $G_2$ the group
$Aut(k)$ of automorphisms of the line, i.e. the set of morphisms of the form
$z\mapsto az+b$, where $a\not=0$. We endow $X_{n,m}$ with two algebraic
actions $\varphi,\psi$ of $G_1,G_2$ respectively, defined by the
formulas:
$$
\begin{array}{rcl}
\varphi_{t}(x_1,..,x_n,y_1,..,y_m,z) &=&(x_1,..,x_n,y_1,..,y_m,z+ ty_1) \\
\psi_{(a,b)}(x_1,..,x_n,y_1,..,y_m,z) & =&(x_1,..,x_n,ay_1,..,ay_m,z+ by_1)
\end{array}
$$
Since the morphisms $\varphi_t ^*$ and $\psi_{(a,b)}^*$ preserve the ring $A_{n,m}$, $\varphi$
and $\psi$ induce two actions of $G_1$ and $G_2$ on the variety $Y_{n,m}$. Moreover, the map
$F$ is equivariant for both $G_1$ and $G_2$. It is then easy to check that:
$$
B_{n,m} ^{G_1}=k[x_1,..,x_n,y_1,..,y_m] \quad \mbox{and} \quad B_{n,m} ^{G_2}=k[x_1,...,x_n]
$$
\begin{thh} \label{second}
Let $X_{n,m}$, $Y_{n,m}$, $\varphi$ and $\psi$ be the varieties and actions defined above. Then the singular set of $Y_{n,m}$ has dimension $\leq n$.
Moreover the actions $\varphi$ and $\psi$ enjoy the following properties:
\begin{itemize}
\item{$A_{n,m}^{G_1}=k[x_1 ^{i_1}..x_n ^{i_n}y_1 ^{j_1}..y_m ^{j_m}, j_1+..+j_m>0]$. In particular, $A_{n,m}^{G_1}$ is not
finitely generated, $B_{n,m} ^{G_1}$ is algebraic but not integral over $A_{n,m}^{G_1}$ and $F^{G_1}$ is
quasifinite but not finite.}
\item{$A_{n,m}^{G_2}=k$. In particular, $B_{n,m} ^{G_2}$ is not algebraic
over $A_{n,m}^{G_2}$, the transcendence degree of $B_{n,m}^{G_2}$ over $A_{n,m} ^{G_2}$ is equal to
$n$ and $F^{G_2}$ is not even quasifinite.}
\end{itemize}
\end{thh}
First note that the variety $Y_{1,m}$ has dimension $m+2$ and its singular set $Sing(Y_{1,m})$ is a line. In particular,
this line consists solely of non normal singular points, and this shows that theorem \ref{premier} is optimal in terms of
singularities. In Van den Essen's example, the second variety has a singular set of dimension 2. Second note that the
difference of dimension between the quotient varieties for the action of $G_2$ is equal to $n$, hence it can be chosen arbitrarily large.

\section{Proof of theorem \ref{premier}}

In this section, we are going to give a proof of theorem \ref{premier}, first in the case when $A$ is a normal ring, and then in the general
case. We begin by recalling an elementary result from Differential Galois Theory, which can be found for instance in \cite{Ko}.

\begin{lem} \label{ext}
Let $k\subset K \subset L$ be fields of characteristic zero, such that the extension $L/K$ is finite. If $d$ is a $k$-derivation on $K$, then
$d$ extends uniquely to a $k$-derivation $D$ of $L$.
\end{lem}

\begin{lem} \label{ext2}
Let $A\subset B$ be two finitely generated integral $k$-algebras, where $char(k)=0$, such that $B$ is integral over $A$. Let
$\FF$ be a family of $k$-derivations on $B$ such that $d(A)\subset A$ for any $d$ in $\FF$. If $A$ is normal, then $B^{\FF}$
is integral over $A^{\FF}$.
\end{lem}
\dem Denote by $K$ the fraction field of $A$, and by $K'$ that of $B$.
Since $B$ is integral over $A$ and that $A$, $B$ are finitely generated $k$-algebras, $B$ is a finite $A$-module. In particular,
$K'/K$ is a finite extension and $k$ is contained in both $K$ and $K'$. Denote by $L/K$
a finite Galois extension containing $K'/K$, with Galois group $G$. Such an extension exists since $K$ has characteristic zero.
For any element $x$ of $B^{\FF}$, consider the polynomial:
$$
P(t)=\prod_{g\in G} (t-g(x))
$$
It is clear that the coefficients of $P$ are invariant with respect to $G$, hence they belong to $K$. Since $x$ is integral over $A$, $g(x)$ is integral
over $A$ for all $g\in G$, and the coefficients of $P$ are integral over $A$. Since they all belong to $K$, and that $A$ is normal, they all lie in $A$.
There remains to show that the coefficients of $P$ are annihilated by every element of $\FF$.

Let $d$ be any derivation belonging to $\FF$. Then $d$ defines a $k$-derivation on $K$. Denote by $D$ its unique extension to $L$ (see lemma \ref{ext}). For
any $g\in G$, consider the map:
$$
D_g=g^{-1}\circ D \circ g
$$
Since $k\subset K$ and $K$ is $G$-invariant, $D_g$ is $k$-linear and $D_g$ coincides with
$d$ on $K$. Moreover $D_g$ is a $k$-derivation on $L$. Indeed, for any $x,y\in L$, we have:
$$
\begin{array}{rcl}
D_g(xy)&=&g^{-1}\circ D(g(x)g(y)) \\
&=&g^{-1}(g(x)D(g(y)) + g(y)D(g(x))) \\
&=&xg^{-1}D(g(y)) + y g^{-1}D(g(x))\\
&=&xD_g(y) + y D_g(x)
\end{array}
$$
By uniqueness of the extension, $D=D_g$ on $L$. In particular, $D\circ g=g\circ D$ for any $g\in G$. Since $D(x)=d(x)=0$, we find $D(g(x))=0$
for all $g\in G$. So the coefficients of $P$ all lie in the kernel of $d$. Since this holds for any derivation in $\FF$, these coefficients
all belong to $A^{\FF}$ and the result follows.
\qed

\begin{lem} \label{ext3}
Let $A$ be a finitely generated integral $k$-algebra. Let $A'$ be its integral closure. If $A$ has isolated
non normal singular points, then $A'/A$ is a finite dimensional $k$-vector space.
\end{lem}
\dem Let $\{x_1,...,x_n\}$ be the collection of non normal singular points of $Spec(A)$, viewed as maximal ideals of $A$,
and set $I=x_1 \cap ...\cap x_n$. First we claim that $A/I$ is finite dimensional over $k$. Indeed since the $x_i$ are
maximal ideals, we obtain by the Chinese Remainder Theorem:
$$
\frac{A}{I} \simeq \frac{A}{x_1} \times ...\times \frac{A}{x_n}
$$
Each $k$-algebra $A/x_i$ is finitely generated and is a field, hence it is a finite extension of $k$ (see \cite{Hu}). In particular,
every quotient $A/x_i$ has finite dimension, and $dim_k A/I <+\infty$.

Second we show that $A/I^m$ is finite dimensional
over $k$ for any $m>0$. The $A$-module $A/I^m$ is filtered by the submodules $M_i=I^i/I^m$ for $i=0,...,m$, and
$M_i/M_{i+1}=I^i/I^{i+1}$. Since $A$ is noetherian and $I$ is an ideal of $A$, $M_i/M_{i+1}$ is a finite $A/I$-module,
hence finite dimensional over $k$. Therefore we have $dim_k A/I^m <+\infty$.

Eventually we prove that $A'/A$ is finite dimensional over $k$. Let $f_1,...,f_r$ be a set of nonzero generators of
$I$. For any $f_i$, the localization
$A_{(1/f_i)}$ has no non normal singular points, hence it is a normal ring and $A\subset A'\subset A_{(1/f_i)}$.
Since $A'$ is a finite $A$-module, there exists
an integer $n_i$ such that $f_i ^{n_i} A'\subset A$. If $m=r\max\{n_i\}$, then $I^m A'\subset A$ and $I^m(A'/A)=0$.
So $A'/A$ is a finite $A/I^m$-module. Since $dim_k A/I^m <+\infty$, the result follows.
\qed
{\it Proof of theorem \ref{premier}}: Let $A\subset B$ be two finitely generated integral $k$-algebras, where $char(k)=0$, such that
$B$ is integral over $A$. Let $\FF$ be a family of $k$-derivations on $B$ such that $d(A)\subset A$ for any $d$ in $\FF$. Assume that $A$ has
isolated non normal singular points. Let us prove that $B^{\FF}$ is integral over $A^{\FF}$.

Let $K$ be the fraction field of $B$, and consider the $k$-subalgebra $B'$ of $K$ generated by $B$ and the integral closure $A'$ of $A$. Since
$A,B$ are finitely generated, $A'$ and $B'$ are also finitely generated. By construction, $B'$ is integral over $A'$. Moreover, each $k$-derivation
$d\in \FF$ extends to a unique $k$-derivation on $A'$ by Seidenberg's Theorem (see \cite{Sei}). Every derivation $d$ is also well-defined
on $K$. Since $d(A')\subset A'$ and $d(B)\subset B$, we have $d(B')\subset B'$ for any $d\in \FF$. Since $A'$ is normal, $(B')^{\FF}$ is integral over $(A')^{\FF}$ by lemma \ref{ext2}. Since $B^{\FF}\subset (B')^{\FF}$, there only remains to show that $(A')^{\FF}$ is integral over $A^{\FF}$.
Let $x$ be any element of $(A')^{\FF}$. Since $A$ has isolated non normal singular points, $A'/A$ is finite dimensional by
lemma \ref{ext3}. In particular, there exist some elements $a_0,...,a_{n-1}\in k$ such that:
$$
P(x)= x^n + a_{n-1}x^{n-1}+ ....+ a_0 \; \in \; A
$$
By construction, $P(x)$ is annihilated by all elements of $\FF$. So $P(x)$ belongs to $A^{\FF}$ and $x$ satisfies
an integral relation with coefficients in $A^{\FF}$. Since this holds for any $x$ in $(A')^{\FF}$,  $(A')^{\FF}$
is integral over $A^{\FF}$ and the result follows.
\qed

\section{A lemma on unipotent group actions}

\begin{lem} \label{uni}
Let $G$ be a unipotent algebraic group over an algebraically closed field $k$. Let $F:X\rightarrow Y$ be a finite morphism
of affine irreducible $G$-varieties. Then $F^G$ is quasifinite.
\end{lem}
\dem Set $B=\OX$ and $A=\OY$. The morphism $F$
induces an inclusion $A\subset B$ such that $B$ is integral over $A$. Choose an element $x$ of $B$ which is $G$-invariant, and let $P(x)=a_n t^n
+...+a_0$ be a polynomial in $A[t]$, of minimal degree $n$ such that $P(x)=0$.  Consider the subset $M$ of $A^{n+1}$ consisting of the
$(n+1)$-tuples $(a_0,...,a_n)$ such that $a_n x^n+...+a_0=0$. By construction, $M$ is an $A$-submodule of $A^{n+1}$, and $A^{n+1}$
is endowed with
the action of $G$ defined by:
$$
g.(a_0,...,a_n)=(g.a_0,...,g.a_n)
$$
Since $x$ is $G$-invariant, $M$ is a rational $G$-submodule of $A^{n+1}$. Now since $G$ is unipotent, there exists a nonzero
element $b=(b_0,...,b_n)$ of $M$
which is $G$-invariant (see \cite{B-S}). In particular, all the $b_i$ are $G$-invariant and we have:
$$
b_n x^n + ...+b_0=0
$$
Note that $b_n$ cannot be equal to zero. Otherwise, all the $b_i$ would be zero by minimality of $n$, a contradiction. So $x$
is algebraic over $A^G$. Since this holds for any $x$ in $B^G$, $B^G$ is algebraic over $A^G$. In particular, the morphism $F^G$
is quasi-finite.
\qed

\section{Properties of the varieties $Y_{n,m}$}

In this section, we are going to establish theorem \ref{second}. We will begin with
a few lemmas concerning the ring $A_{n,m}$ defined in the introduction.

\subsection{A few preliminary lemmas}

\begin{lem}
$k[x_1 ^{i_1}...x_n^{i_n}y_1 ^{j_1}...y_m ^{j_m}, j_1+...+j_m >0] \subset A_{n,m}$.
\end{lem}
\dem We are going to prove by induction on $p=i_1+...+i_n$ that every monomial of the form $x_1 ^{i_1}...x_n^{i_n}y_1 ^{j_1}...y_m ^{j_m}$, where
$j_1+...+j_m >0$, belongs to $A_{n,m}$. Since $y_1,...,y_n$ belong to $A_{n,m}$, we may restrict ourselves to the monomials of the
form $x_1 ^{i_1}...x_n^{i_n}y_j$, where $j=1,...,m$. For $p=0$, this is clear because $A_{n,m}$ contains $y_1,...,y_m$. Assume the
property holds to the order $p$. For convenience, we set $t_i=x_i ^2 + x_i z$ and note that every $t_i$ belongs to
$A_{n,m}$. Consider any monomial $a$ of the form $x_1 ^{i_1}...x_n^{i_n}y_j$, where
$i_1+...+i_n=p+1$. If all the $i_k$ are $\leq 1$, then $a$ belongs to $A_{n,m}$ by construction. If one of the $i_k$
is $\geq 2$, for instance $i_1\geq 2$, then write $a=x_1 ^2 b$, where $b=x_1 ^{i_1 -2}...x_n^{i_n}y_j$
and set $c=x_1 ^{i_1 -1}...x_n^{i_n}y_j$. Since $x_1 ^2 + x_1 z -t_1=0$, we obtain by multiplication by $b$:
$$
a + zc - t_1 b=0
$$
By our induction's hypothesis, the monomials $b$ and $c$ belong to $A_{n,m}$. So $a$ belongs to $A_{n,m}$
and the result follows.
\qed

\begin{lem} \label{infini}
$k[x_1,...,x_n,y_1,...,y_m]\cap A_{n,m}=k[x_1 ^{i_1}...x_n^{i_n}y_1 ^{j_1}...y_m ^{j_m}, j_1+...+j_m >0]$.
\end{lem}
\dem Let $f$ be an element of $k[x_1,...,x_n,y_1,...,y_m]\cap A_{n,m}$. By the previous lemma, we know that
every monomial of $k[x_1,...,x_n,y_1,...,y_m]$ containing at least one of the $y_j$ belongs to $A_{n,m}$. Up to
substracting a linear combination of such monomials to $f$, we may assume that $f$ only depends on $x_1,...,x_n$.
If we show that such an $f$ is always constant, then the lemma will follow. So we are going to prove by
induction on $n\geq 1$
that any element $f$ of $A_{n,m}$ that only depends on $x_1,...,x_n$ is a constant. For $n=1$, consider such an element
$f=f(x_1)$ of $A_{1,m}$. Then there exists a polynomial $P$ such that:
$$
f(x_1)=P(y_1,...,y_m,z,x_1 ^2 +zx_1,x_1^3 +zx_1^2,x_1y_1,...,x_1y_m) 
$$
By setting $y_1=...=y_m=0$, we can see there exists a polynomial $Q$ such that:
$$
f(x_1)=Q(z,x_1 ^2 +zx_1,x_1^3 +zx_1^2)
$$
If $x_1=0$, then $Q(z,0,0)=f(0)$ is a constant. So the polynomial $Q(a,b,c)-f(0)$ has no pure
terms in $a$, and it can be expanded as:
$$
Q(a,b,c)-f(0)=\sum_{k+l>0} Q_{k,l}(a)b^k c^l
$$
In particular, this yields for $f(x_1)-f(0)$:
$$
f(x_1)-f(0)=\sum_{k+l>0} Q_{k,l}(z)(x_1 ^2 +zx_1)^k(x_1^3 +zx_1^2)^l
$$
Since $x_1 +z$ divides both $x_1 ^2 +zx_1$ and $x_1^3 +zx_1^2$, $x_1 +z$ must divide $f(x_1)-f(0)$, which
is impossible unless $f$ is constant. Now assume the property holds to the order $(n-1)$, and let $f=f(x_1,...,x_n)$ be
an element of $A_{n,m}$. Then there exists a polynomial $P$ such that:
$$
f(x_1,...,x_n)=P(y_1,...,y_m,z,x_1 ^2 +zx_1,x_1^3 +zx_1^2,...,x_n ^2 +zx_n,x_n^3 +zx_n^2,...,x_1 ^{i_1}...x_n ^{i_n}y_j,...) 
$$
By setting $y_1=...=y_m=0$, we can see there exists a polynomial $Q$ such that:
$$
f(x_1,...,x_n)=Q(z,x_1 ^2 +zx_1,x_1^3 +zx_1^2,...,x_n ^2 +zx_n,x_n^3 +zx_n^2)
$$
By setting $x_n=0$, we find:
$$
f(x_1,...,x_{n-1},0)=Q(z,x_1 ^2 +zx_1,x_1^3 +zx_1^2,...,x_{n-1} ^2 +zx_{n-1},x_{n-1}^3 +zx_{n-1}^2,0,0)
$$
So the polynomial $f(x_1,...,x_{n-1},0)$ belongs to $A_{n-1,m}$. By our induction's hypothesis, $f(x_1,...,x_{n-1},0)$ is constant, and
we may assume that:
$$
Q(a,b_1,c_1,...,b_{n-1},c_{n-1},0,0)-f(0,...,0)=0
$$
In particular,
$Q(a,b_1,c_1,...,b_n,c_n)-f(0,...,0)$ has no pure terms in $a,b_1,c_1,...,b_{n-1},c_{n-1}$, and
it can be expanded as:
$$
Q(a,b_1,c_1,...,b_n,c_n)-f(0,...,0)=\sum_{k+l>0} Q_{k,l}(a,b_1,c_1,...,b_{n-1},c_{n-1})b_n ^k c_n ^l
$$
If $q_{k,l}=Q_{k,l}(z,x_1 ^2 +zx_1,x_1^3 +zx_1^2,...,x_{n-1} ^2 +zx_{n-1},x_{n-1}^3 +zx_{n-1}^2)$, then
this yields:
$$
f(x_1,...,x_n)-f(0,...,0)=\sum_{k+l>0} q_{k,l}(x_n ^2 +zx_n)^k(x_n^3 +zx_n^2)^l
$$
Since $x_n +z$ divides both $x_n ^2 +zx_n$ and $x_n^3 +zx_n^2$, $x_n +z$ must divide $f(x_1,...,x_n)-f(0,...,0)$, which
is impossible unless $f$ is constant, and the result follows. 
\qed

\begin{lem} \label{infini2}
The $k$-algebra ${\cal{A}}_{n,m}=k[x_1 ^{i_1}...x_n^{i_n}y_1 ^{j_1}...y_m ^{j_m}, j_1+...+j_m >0]$ is not finitely generated.
\end{lem}
\dem Suppose on the contrary that ${\cal{A}}_{n,m}$ is finitely generated, and let $f_1,...,f_r$ be a system of generators.
For convenience, we may assume that $f_i(0,...,0)=0$ for any $i$. Since every $f_i$ is a linear combination of monomials
of the form $x_1 ^{i_1}...x_n^{i_n}y_1 ^{j_1}...y_m ^{j_m}$, where $j_1+...+j_m >0$, we may even assume that $\{f_1,...,f_r\}$
consists solely of such monomials. Then for any couple of integers $i_1,j_1$, where $j_1>0$, there exists a polynomial
$P$ such that:
$$
x_1 ^{i_1} y_1 ^{j_1}=P(f_1,...,f_r)
$$
Let $g_i$ be the monomial $f_i(x_1,y_1,0,...,0)$ for any $i$. Then, it is easy to check that $g_1,...,g_r$ belong to ${\cal{A}}_{1,1}$.
By setting $x_2=...=x_n=y_2=...=y_m=0$, we find:
$$
x_1 ^{i_1} y_1 ^{j_1}=P(g_1,...,g_r)
$$
In particular, $g_1,...,g_r$ span the $k$-algebra ${\cal{A}}_{1,1}$. Set $g_i = x_1 ^{n_i}y_1 ^{m_i}$ for any $i$, and consider the monomial
$x_1 ^s y_1$, where $s>\max\{n_i\}$. Since $x_1 ^s y_1$ belongs to ${\cal{A}}_{1,1}$, there exists a polynomial $P$ such that:
$$
x_1 ^{s} y_1 =P(g_1,...,g_r)
$$
Since the $g_i$ are monomials, $P$ must contain a monomial of the form $u_1 ^{a_1}...u_r ^{a_r}$ such that:
$$
\begin{array}{rcl}
s&=&a_1 n_1 +...+a_r n_r \\
1&=&a_1 m_1+...+a_r m_r
\end{array}
$$
Since the $a_i$ are nonnegative integers and that $m_i>0$ for any $i$, this implies that $a_i$ is zero for every index $i$ except
one, say $i_0$, and that $a_{i_0}=1$. But then $s=n_{i_0}$, which is impossible since $s>n_i$ for any $i$.
\qed

\subsection{Proof of theorem \ref{second}}

Let $Y_{n,m}$ be the variety defined in the introduction. We first show that $Sing(Y_{n,m})$ has dimension $\leq n$.
If we localize $A_{n,m}$ with respect to any $y_i$, then we get the ring:
$$
(A_{n,m})_{\frac{1}{y_i}}=k[x_1,...,x_n,y_1,..,y_i, \frac{1}{y_i},..,y_m, z]
$$
which is obviously smooth. Again by localizing $A_{n,m}$ with respect to either $x_1 ^{i_1}...x_n ^{i_n}y_j$, or $(x_1 ^2 + x_1 z)...(x_n ^2 + x_n z)$,
we get a regular ring. This implies that $Sing(Y_{n,m})$ is contained in the zero set of the ideal $I$
generated by all these polynomials. In particular, $Sing(Y_{n,m})$ is contained in the image by $F$
of the zero set $E$ of all these polynomials in $k^{n+m+1}$. This set $E$ consists of a finite union
of sets of the form $V(y_1,...,y_n,x_i)$ and $V(y_1,...,y_n,x_i+z)$. So $E$ has dimension $n$
and we get:
$$
dim\; Sing(X_{n,m})\leq n
$$
Second we compute the ring of invariants $A_{n,m} ^{G_1}$. By construction, we clearly have $B_{n,m} ^{G_1}=k[x_1,..,x_n,y_1,..,y_m]$. Since $A_{n,m}$ is a subalgebra of $B_{n,m}$, we obtain by lemma \ref{infini}:
$$
A_{n,m} ^{G_1}= A_{n,m} \cap k[x_1,..,x_n,y_1,..,y_m]= k[x_1 ^{i_1}...x_n^{i_n}y_1 ^{j_1}...y_m ^{j_m}, j_1+...+j_m >0]
$$
By lemma \ref{infini2}, this algebra is not finitely generated. We claim that $B_{n,m} ^{G_1}$ cannot be integral
over $A_{n,m} ^{G_1}$. Indeed, assume that $B_{n,m} ^{G_1}$ is integral over $A_{n,m} ^{G_1}$. Then the
$x_i$, $y_j$ satisfy some integral relations over $A_{n,m} ^{G_1}$. Let $A_0$ be the $k$-algebra generated
by all the coefficients of these relations. By construction, $A_0$ is a finitely generated subalgebra of $A_{n,m} ^{G_1}$, and
$B_{n,m} ^{G_1}$ is integral over $A_0$. Since $B_{n,m} ^{G_1}$ is finitely generated, $B_{n,m} ^{G_1}$ is a finite
$A_0$-module. But $A_0$ is noetherian, and $A_{n,m} ^{G_1}$ is an $A_0$-submodule of $B_{n,m} ^{G_1}$, so $A_{n,m} ^{G_1}$
is a finite $A_0$-module. If $\{e_1,...,e_r\}$ denotes a basis of $A_{n,m} ^{G_1}$ over $A_0$, then:
$$
A_{n,m} ^{G_1}=A_0[e_1,...,e_r]
$$
In particular, $A_{n,m} ^{G_1}$ is finitely generated, hence a contradiction. However, $B_{n,m} ^{G_1}$ is algebraic over $A_{n,m} ^{G_1}$. Indeed,
they have the same fraction field $k(x_1,..,x_n,y_1,..,y_m)$. 

Consider now the ring $B_{n,m} ^{G_2}$. The group $G_2$ is the semi-direct product of $G_a(k)$ and $G_m(k)$. The action of $G_a(k)$ corresponds to
the previous action $\varphi$ of $G_1$. The action of $G_m(k)$ is related to the weighted homogeneous degree $deg$ on $k[x_1,..,x_n,y_1,..,y_m,z]$,
which assigns the weight $0$ to each $x_i$ and $z$, and $1$ to each $y_i$. In particular, the invariants of $G_2$ on $k[x_1,..,x_n,y_1,..,y_m,z]$
are the $G_1$-invariant polynomials of degree zero with respect to the $y_i$. More precisely:
$$
B_{n,m} ^{G_2}=k[x_1,...,x_n]
$$
By lemma \ref{infini}, the $G_1$-invariants of $A_{n,m}$ of degree zero with respect to the $y_i$ are the constants, i.e.
$A_{n,m} ^{G_2}=k$. In particular, $B_{n,m} ^{G_2}$ is not even algebraic over $A_{n,m} ^{G_2}$, and the transcendence degree of $B_{n,m} ^{G_2}$ over $A_{n,m} ^{G_2}$
is equal to $n$. This ends the proof of theorem \ref{second}.


\begin{thebibliography}{D-F}

\bibitem[B-S]{B-S} M.Brion, G.Schwarz {\it Th\'eorie des invariants et G\'eom\'etrie des vari\'et\'es quotients}, Collection Travaux en Cours 61, Hermann, Paris, 2000.

\bibitem[D-F]{D-F} D.Daigle, G.Freudenburg {\it A counterexample to Hilbert's fourteenth problem in dimension $5$},
J. Algebra 221 (1999), $n^o 2$, 528-535.

\bibitem[Ei]{Ei} D.Eisenbud {\it Commutative Algebra with a view toward
Algebraic Geometry}, Graduate Texts in Mathematics, 150, Springer Verlag New York, 1995.

\bibitem[Hu]{Hu} J.Humphreys {\it Linear algebraic groups}, Graduate Texts in Mathematics $n^o$ 21, Springer Verlag New-York Heidelberg, 1975.

\bibitem[Ko]{Ko} E.R.Kolchin {\it Differential Algebra and Algebraic Groups}, Pure and Applied Mathematics 54,
Academic Press, New York London, 1973.

\bibitem[Kr]{Kr} H.Kraft {\it Geometrische Methoden in der Invariantentheorie}, Aspects of Mathematics, D1 Friedr.
Vieweg $\&$ Sohn, Braunschweig, 1984.

\bibitem[Na]{Na} M.Nagata {\it On the fourteenth problem of Hilbert}, 1960 Proc. Internat. Congress Math. 1958, pp. 459-462,
Cambridge Univ.Press, New York.

\bibitem[Sei]{Sei} A.Seidenberg {\it Derivations and integral closure}, Pacific J. Maths
16, 1966, 167-173.

\bibitem[Sh]{Sh} I.Shafarevich {\it Basic Algebraic Geometry 1: Varieties in projective spaces}, second edition,
Springer Verlag, Berlin, 1994.

\bibitem[Van]{Van} A.Van den Essen {\it Polynomial Automorphisms and the Jacobian Conjecture}, Progress in Maths 190,
Birkh\"auser Verlag, Basel 2000.

\end{thebibliography}
\end{document}